
\font\eusm=eusm10                   


\font\eusms=eusm7                       

\font\eusmss=eusm5                      

\font\scriptsize=cmr7

\input amstex

\documentstyle{amsppt}
  \magnification=1000
  \hsize=7.0truein
  \vsize=9.0truein
  \hoffset -0.1truein
  \parindent=2em

\define\Ac{{\Cal A}}                         

\define\Ad{\text{\rm Ad}}                    

\define\Afr{{\frak A}}                       

\define\Afro{\overset{\scriptsize o}\to\Afr} 

\define\ah{\hat a}                           

\define\alphah{\hat\alpha}                   

\define\Aut{\text{\rm Aut}}                  

\define\Bc{{\Cal B}}                         

\define\bh{\hat b}                           

\define\bigoplustxt{{                        
     \tsize\bigoplus}}

\define\bigotimestxt{{                       
     \tsize\bigotimes}}

\define\Bof{{\eusm B}}                       

\define\Cc{{\Cal C}}                         

\define\Cpx{\bold C}                         

\define\dif{\text{\it d}}                    

\define\ef{\text{\rm ef}}                    

\define\eh{\hat e}                           

\define\freeF{\bold F}                       

\define\freeprodi{\mathchoice                
     {\operatornamewithlimits{\ast}
      _{\iota\in I}}
     {\raise.5ex\hbox{$\dsize\operatornamewithlimits{\ast}
      _{\sssize\iota\in I}$}\,}
     {\text{oops!}}{\text{oops!}}}

\define\Gh{\hat G}                           

\define\Hil{\mathchoice                      
     {\text{\eusm H}}
     {\text{\eusm H}}
     {\text{\eusms H}}
     {\text{\eusmss H}}}

\define\HilA{\Hil_\Ac}                       

\define\HilAbar{\Hilbar_\Ac}                 

\define\Hilbar{\overline\Hil}                

\define\HilN{\Hil_\NvN}                      

\define\HilNbar{\Hilbar_\NvN}                

\define\HilNo{\Hilo_\NvN}                    

\define\Hilo{{\overset{\scriptsize o}        
     \to\Hil}}

\define\id{\text{\it id}}                    

\define\Int{\text{\rm Int}}                  

\define\Integers{\bold Z}                    

\define\MvN{{\Cal M}}                        

\define\MvNo{\overset{\scriptsize o}\to\MvN} 

\define\ndn{                                 
     \neq \cdots \neq }

\define\nm#1{||#1||}                         

\define\NvN{{\Cal N}}                        

\define\otdt{\otimes\cdots\otimes}           

\define\pibar{\overline\pi}                   

\define\QED{$\hfill$\qed\enddemo}            

\define\Real{{\bold R}}                      

\define\restrict{\lower .3ex                 
     \hbox{\text{$|$}}}


\define\tdt{                                 
     \bigotimestxt\cdots\bigotimestxt}

\define\Tr{\text{\rm Tr}}                    

\define\uh{\hat u}                           

\define\xh{\hat x}                           

\define\xIm{x_{\text{\rm Im}}}               

\define\xRe{x_{\text{\rm Re}}}               




\topmatter

  \title Factoriality and the Connes invariant $T(\MvN)$
         for free products of von Neumann algebras \endtitle

  \author Kenneth J. Dykema \endauthor

  \thanks This work was supported by the
          Fannie and John Hertz Foundation and a National
          Science Foundation Postdoctoral Fellowship. \endthanks

  \affil Department of Mathematics \\ University of California \\
         Berkeley, CA 94720 \endaffil

  \date Revised Version, $\;$ July 27, 1993 \enddate

  \abstract Sufficient conditions for factoriality are given for free products
    of von Neumann algebras with respect to states that are not necessarily
    traces.
    The Connes $T$--invariant of the free product algebra is found, which has
    implications for the type of the algebra.
    Roughly speaking, the free product of von Neumann algebras
    with respect to states, one of which is not a trace, is a type~III factor.
   \endabstract

\endtopmatter

\document
  \TagsOnRight
  \baselineskip=14pt

\noindent{\bf Introduction.}
  In this paper, we examine free products of (not necessarily injective)
von Neumann algebras with respect to states that, in contrast to the cases
studied in~\cite{D3}, need not be traces.
Our main theorem gives sufficient conditions for the free product to
be a factor, and using it we are also able to make conclusions
about the type of
the free product factor.
Similarly to the results of~\cite{D3}, for factoriality
of the free product it suffices that
the original algebras not be too ``lumpy,''
{\it i.e\.} not have big minimal projections.
However, in contrast to~\cite{D3}, the sufficient conditions given here
are far from
necessary for factoriality.
Connes' $T$-invariant is computed for the free product factor in terms
of the original algebras and states, and this allows one to determine,
in the case of separable preduals, if the free product factor is
type~III, and to gain information about its classification
as type III$_\lambda$.
Roughly speaking, the pattern for
the type classification obtained is that as long as the states with
respect to which
we take the free product are not all traces, we get a type~III
factor, and it is type III$_\lambda$, where we are able to specify
$\lambda$ to within a choice of two values, one of which is $0$.

  Since the appearance of Voiculescu's
free probability theory in the mid~'80's,
(see~\cite{VDN} and references contained therein),
and more specifically since Voiculescu's random matrix model for
freeness~\cite{V3}, there has been steady progress in understanding
the free group factors and free products of finite von Neumann
algebras, \linebreak \cite{V2,R1,R2,R3,R4,D1,D2,D3}.
For example, in~\cite{D3} the free product of any two finite injective
von Neumann algebras was found and expressed in terms of the interpolated
free group factors~\cite{D2}, (see also~\cite{R3}), via
``free dimension.''

  More recently, in~\cite{R4} F\. R\u{a}dulescu exhibited a 1--parameter
group of trace--scaling automorphisms of $L(\freeF_\infty)\otimes\Bof(\Hil)$.
In~\cite{R5}, he showed that if $\psi_\lambda$ is the state on the $2\times2$
matrices $M_2(\Cpx)$ defined with respect to a system of matrix units
$\{e_{\iota j}\mid\iota,j=1,2\}$ by $\psi_\lambda(e_{12})=0$,
$\psi_\lambda(e_{11})=\frac1{1+\lambda}$ and
$\psi_\lambda(e_{22})=\frac\lambda{1+\lambda}$
for some $0<\lambda<1$, then the free product of a diffuse abelian von Neumann
algebra with $(M_2(\Cpx),\psi_\lambda)$
is a type~III$_\lambda$ factor having a discrete decomposition
given by the appropriate subgroup of the 1--parameter group from~\cite{R4}.

  Note also the related work of L\. Barnett~\cite{B}.

\medskip

  One of A\. Connes' great achievements is his classification of type~III
factors in terms of trace scaling automorphisms of type II$_\infty$
factors~\cite{C}.
The type~III factors are differentiated as type III$_\lambda$,
($0\le\lambda\le1$).
Recall that his invariant $S(\MvN)$ is the intersection of the spectra
of the modular operators
$\Delta_\phi$, as $\phi$ ranges over the set of all normal,
semifinite, faithful (n.s.f\.) weights on $\MvN$.
The Connes classification is then that
$$ \MvN\text{ is type }\cases\text{ I or II}\quad&\text{if }S(\MvN)=\{1\} \\
\text{ III}_0&\text{if }S(\MvN)=\{0,1\} \\
\text{ III}_\lambda&\text{if }S(\MvN)
=\{0\}\cup\{\lambda^n\mid n\in\Integers\},\;
0<\lambda<1 \\
\text{ III}_1&\text{if }S(\MvN)=[0,+\infty).
\endcases $$

  Connes defined another invariant
$T(\MvN)=\{t\in\Real\mid\sigma^\phi_t\in\Int\MvN\}$,
where $(\sigma^\phi_t)_{t\in\Real}$ is the modular automorphism
group of $\MvN$ corresponding to the
n.s.f\. weight
$\phi$, and $\Int\MvN$ is the group of inner automorphisms of $\MvN$.
He proved that this set is independent of $\phi$
and also that if $t\in T(\MvN)$, then there is an n.s.f\. weight
$\psi$ such that $\sigma^\psi_t=\id$.
This invariant
$T(\MvN)$ is related to the type classification of factors as follows:
$$ T(\MvN)=\cases\Real&\text{if $\MvN$ is type I or II} \\
\{0\}&\text{if $\MvN$ is type III}_1 \\
\frac{2\pi}{\ln\lambda}\Integers\quad&\text{if $\MvN$ is type III}_\lambda,\;
\lambda\in(0,1).
\endcases $$
If a von Neumann algebra $\MvN$
has separable predual, then $T(\MvN)=\Real$ if and only if $\MvN$ is
semifinite, (see \cite{S}, 27.2).
Thus for a given factor $\MvN$ with separable predual,
knowing $T(\MvN)$ allows one to decide
if $\MvN$ is type III, and if so to specify the $\lambda$--classification to
within a choice of at most two, one of which is III$_0$.

\proclaim{Definition 0}\rm
Let $\phi$ be an n.s.f\. weight on
a von Neumann algebra $\MvN$, and let $(\sigma^\phi_t)_{t\in\Real}$
be the resulting modular automorphism group of $\MvN$.
We define
$$ I(\MvN,\phi)=\{t\in\Real\mid\sigma^\phi_t=\id\}, $$
which is clearly an additive subgroup of $T(\MvN)$.
Note that since $t\mapsto\sigma_t^\phi$ is continuous in the
pointwise--strong topology on $\Aut(\MvN)$, we have that $I(\MvN,\phi)$
is a closed subgroup of $\Real$.
\endproclaim

  Our main theorem shows that for ``most'' families of von Neumann
algebras with states $((\MvN_\iota,\phi_\iota))_{\iota\in I}$,
if $(\MvN,\phi)=\freeprodi(\MvN_\iota,\phi_\iota)$ is their free product,
then
$T(\MvN)=I(\MvN,\phi)=\bigcap_{\iota\in I}I(\MvN_\iota,\phi_\iota)$.

  In \S1, we prove a certainly well known result about the modular
automorphism group of a free product.
In~\S2 we define a technical quantity, called the expansion factor
of a von Neumann algebra with a specified state,
prove some lemmas about this quantity
and prove our main theorem.
In~\S3 the expansion factors of
several von Neumann algebras are calculated,
allowing one to decide
when the hypotheses of our main
theorem are satisfied.
In~\S4 appear examples of specific results
that follow from the main theorem, and also some limitations of
the main theorem.

\noindent{\bf \S1. The modular automorphism group of a free product}

  Let us review the construction
of the free product of von Neumann algebras~\cite{V1}
and then show that the modular automorphism
(with respect to the free product state)
is the free product of the modular automorphisms of the original algebras.
This result has surely been known since Voiculescu in~\cite{V1} found
the modular operator for a free product algebra.
Let $(\MvN_\iota,\phi_\iota)_{\iota\in I}$ be a family of von Neumann
algebras with specified faithful normal states.
Our notation for the GNS construction is
$\Hil_\iota=L^2(\MvN_\iota,\phi_\iota)$ with
$\MvN_\iota\ni x \mapsto\xh\in\Hil_\iota$, and the
distinguished vector $\hat1$ is denoted by $\xi_\iota$.
The Hilbert algebra $\Afr_\iota=(\MvN_\iota)\hat{\,}=\MvN_\iota\xi_\iota$
is dense in $\Hil_\iota$, and we have the left representation
$\pi_\iota$
of $\MvN_\iota$ on $\Hil_\iota$
given by $\pi_\iota(a)\xh=(ax)\hat{\,}$.
The modular operator, $\Delta_\iota$ on $\Hil_\iota$, is the unbounded
positive operator constructed in the following way (see~\cite{T,SZ}).
The conjugation operator,
$S_\iota$, is the closure of the unbounded operator defined by
$S_\iota(\xh)=(x^*)\hat{\,}$.
Then $\Delta_\iota=S_\iota^*S_\iota$.
Note that $\Delta_\iota\xi_\iota=\xi_\iota$.

  Consider the free product of Hilbert spaces (with distinguished
vectors)
$(\Hil,\xi)=\freeprodi(\Hil_{\iota},\xi_{\iota})$,
given by
$$ \Hil = \Cpx\xi \; \bigoplustxt \; \bigoplus_{n \geq 1}
(\bigoplus_{(\iota_{1} \neq \iota_{2} \ndn \iota_{n})}
\Hilo_{\iota_{1}} {\tdt}
\Hilo_{\iota_{n}} ), \tag1 $$
where $\Hilo_{\iota}$ is the orthocomplement
of $\xi_{\iota}$ in $\Hil_{\iota}$.
We often identify $\Hil_\iota$ with $\Cpx\xi\oplus\Hilo_\iota\subseteq\Hil$.
There are injective, unital, normal $*$--homomorphisms
$\lambda_\iota:\MvN_\iota\rightarrow\Bof(\Hil)$, where
$\lambda_\iota(a)$ acts essentially ``on the first component''
of each tensor product in~(1).
To be precise,
for $x\in\ker\phi_\iota$ and
$\zeta_1\otdt\zeta_n\in\Hilo_{\iota_{1}} {\tdt} \Hilo_{\iota_{n}}$
we have
$$ \aligned \lambda_\iota(x)\xi&=\xh \\
\lambda_\iota(x)(\zeta_1\otdt\zeta_n)&=
\cases c\zeta_2\otdt\zeta_n + (\pi_\iota(x)\zeta_1-c\xi_\iota)
\otimes\zeta_2\otdt\zeta_n
&\text{ if }\iota_1=\iota, \text{ where }
c=\langle \pi_\iota(x)\zeta_1,\xi_1\rangle , \\
\xh\otimes\zeta_1\otdt\zeta_n &\text{ if }\iota_1\neq\iota,
\endcases \endaligned \tag2  $$
and of course in the first case of
the second equation above, if $n=1$ we take
$c\xi + (\pi_\iota(x)\zeta_1-c\xi_\iota)$.
Then the free product von Neumann algebra is
$\MvN=(\bigcup_{\iota\in I}\lambda_\iota(\MvN_\iota))''$ and has
free product state $\phi(x)=\langle x\xi,\xi\rangle $.
We write
$(\MvN,\phi)=\freeprodi(\MvN_\iota,\phi_\iota)$.
For ease of notation we often suppress the $\lambda_\iota$ and
simply think of $\MvN_\iota\subset\MvN$.
We sometimes write $x=\pi(x)$ for $x\in\MvN$ acting on $\Hil$, to emphasize the
fact that $\MvN$ acts ``on the left.''
Then $\Hil=L^2(\MvN,\phi)$ and $\Afr=\MvN\Hat{\;\;}=\MvN\xi$ is a Hilbert
algebra containing the dense subspace
$$ \Cpx\xi \; \bigoplustxt \; \bigoplus_{n \geq 1}
(\bigoplus_{(\iota_{1} \neq \iota_{2} \ndn \iota_{n})}
\Afro_{\iota_{1}} {\tdt}
\Afro_{\iota_{n}} ), $$
where $\Afro_\iota=\Afr_\iota\bigcap\Hilo_\iota$.
Now
$$ \Afro_{\iota_{1}} {\tdt} \Afro_{\iota_{n}}
= \MvNo_{\iota_1}\cdots\MvNo_{\iota_n}\xi, $$
where $\MvNo_{\iota_k}=\ker\phi_{\iota_k}$,
so the conjugation operator, $S$, for $\Afr$ acts on
$\zeta_1\otdt\zeta_n\in\Afro_{\iota_{1}} {\tdt} \Afro_{\iota_{n}}$
by
$$ S(\zeta_1\otdt\zeta_n)=(S_{\iota_n}\zeta_n)\otdt(S_{\iota_1}\zeta_1). $$
Thus
$$ S^*(\zeta_1\otdt\zeta_n)
=(S^*_{\iota_n}\zeta_n)\otdt(S^*_{\iota_1}\zeta_1) $$
and so
$$ \Delta(\zeta_1\otdt\zeta_n)=
(\Delta_{\iota_1}\zeta_1)\otdt(\Delta_{\iota_n}\zeta_n), \tag3 $$
(this was first observed in~\cite{V1}).

\proclaim{Theorem 1}
Let $(\MvN,\phi)=\freeprodi(\MvN_\iota,\phi_\iota)$ be as above and
let us write $\MvN_\iota\subset\MvN$.
Let $(\sigma^{\phi_\iota}_t)_{t\in\Real}$ be the modular automorphism
group of $\MvN_\iota$ with respect to the state $\phi_\iota$
and let $(\sigma_t)_{t\in\Real}$ be the modular automorphism
group of $\MvN$ with respect to the free product state $\phi$.
Then $\sigma_t(x)=\sigma^{\phi_\iota}_t(x)$
$\forall x\in\MvN_\iota$, which of course determines $\sigma_t$.
\endproclaim
\demo{Proof}
The modular automorphism is defined by $\sigma_t(x)=\Delta^{-it}x\Delta^{it}$
$\forall x\in\MvN$
and $\sigma^{\phi_\iota}_t(x)=\Delta_\iota^{-it}x\Delta_\iota^{it}$
$\forall x\in\MvN_\iota$.
Writing now the injections $\lambda_\iota:\MvN_\iota\rightarrow\MvN$
to be absolutely clear, we must show that
$\sigma_t(\lambda_\iota(x))=\lambda_\iota(\sigma^{\phi_\iota}_t(x))$
$\forall x\in\MvN_\iota$.
Since every modular automorphism
sends $1$ to $1$, we can concentrate on the
case $x\in\MvNo_\iota$.
But using~(2) and~(3) gives $\Delta^{-it}\lambda_\iota(x)\Delta^{it}\xi
=\Delta^{-it}\lambda_\iota(x)\xi=\Delta^{-it}\xh=\Delta_\iota^{-it}\xh
=\Delta_\iota^{-it}x\xi_\iota=\Delta_\iota^{-it}x\Delta_\iota^{it}\xi_\iota
=\sigma^{\phi_\iota}_t(x)\xi_\iota=(\sigma^{\phi_\iota}_t(x))\hat{\,}
=\lambda_\iota(\sigma^{\phi_\iota}_t(x))\xi.$
But two elements of $\MvN$ that agree on $\xi$ must be equal.
\QED

\noindent{\bf \S2 Factoriality and the $T$--invariant of a free product.}

  Let $(\MvN,\phi)$ be a von Neumann algebra with normal, faithful (n.f\.)
state,
let $\Hil=L^2(\MvN,\phi)$ have distinguished vector $\xi=\hat1$
and set $\Hilo=\Hil\ominus\Cpx\xi$.
We will denote by $\Hilbar$ the dual Hilbert space of $\Hil$, which
is just $\Hil$ with conjugate scalar multiplication.
$P_\xi$ will denote the orthogonal projection of $\Hil$ onto $\Cpx\xi$
and $P_\Hilo=I-P_\xi$ the orthogonal projection of $\Hil$ onto $\Hilo$.
$\MvN$ acts by bounded operators on the left of $\Hil$ by
$\pi(a)\xh=(ax)\hat{\,}$ $\forall a,x\in\MvN$.
Hence there is the dual anti--representation $\pibar$ of $\MvN$ on $\Hilbar$
given by
$$ \langle\zeta,\pibar(x)\eta\rangle=\langle\pi(x)\zeta,\eta\rangle
\;\forall\,x\in\MvN,\zeta\in\Hil,\eta\in\Hilbar. $$
Of course, viewing elements of $\Hilbar$ as being elements of $\Hil$,
we have $\pibar(x)\eta=\pi(x^*)\eta$.

\proclaim{Definition 2.1}\rm
Recall $\Hil\bigotimes\Hil$ is equal to the space of Hilbert--Schmidt
operators from $\Hilbar$ into $\Hil$.
Consider a bounded linear map $\alpha:\Hilbar\rightarrow\MvN$.
For $v\in\Hil$ we thus have
$$ \pi(\alpha(\cdot))v:\Hilbar\rightarrow\Hil. \tag4 $$
Define $\Hil\otimes\MvN$ to be the space of $\alpha$
such that the operator in~(4) is Hilbert--Schmidt,
{\it i.e\.} $\pi(\alpha(\cdot))v\in\Hil\otimes\Hil$,
for every $v\in\Hil$, and furthermore that the norm of
$\pi(\alpha(\cdot))v$ in $\Hil\otimes\Hil$, ({\it i.e\.} its
Hilbert--Schmidt norm), is uniformly bounded for $\nm{v}=1$.
We then say that $\nm{\alpha}_{HS}=
\sup_{\nm{v}=1}\nm{\pi(\alpha(\cdot))v}_{HS}$.
Set $\alphah=\pi(\alpha(\cdot))\xi\in\Hil\otimes\Hil$.
For $\zeta\in\Hil$ and $x\in\MvN$ the simple tensor
$\zeta\otimes x\in\Hil\otimes\MvN$ is given by
$(\zeta\otimes x)(\eta)=\langle\zeta,\eta\rangle x$.
Note that for $v\in\Hil$ we have $\pi((\zeta\otimes x)(\cdot))v
=\zeta\otimes(\pi(x)v)\in\Hil\otimes\Hil$.

  A bounded linear map $\alpha:\Hilbar\rightarrow\MvN$ is said to be
{\it left--right $\MvN$--equivariant} if
$$ \alpha(\pibar(a)\eta)=\alpha(\eta)a\;\;\forall\,a\in\MvN,\eta\in\Hilbar. $$
The {\it expansion factor} of $(\MvN,\phi)$, denoted $\ef(\MvN,\phi)$,
is defined to be the largest constant $c\ge0$ such that
$$ \nm{(P_{\Hilo}\otimes P_{\Hilo})\alphah}
\ge c|\langle\alphah,\xi\otimes\xi\rangle| \tag5 $$
for every left--right $\MvN$--equivariant $\alpha\in\Hil\otimes\MvN$.
Of course, if the right hand side of~(5) is always zero,
we have $\ef(\MvN,\phi)=+\infty$.
\endproclaim

\proclaim{Remark 2.2} \rm
The definitions above become somewhat more transparent,
and the expansion factor easier to calculate, in the case
where $\MvN$ is finite dimensional.  For then every linear
map $\alpha:\Hilbar\rightarrow\MvN$ is a sum of simple
tensors, hence $\in\Hil\otimes\MvN$.  The definition of left--right
$\MvN$--equivariance can be rephrased in terms of $\alphah$ only.
For $a\in\MvN$, let $\rho(a)$ act ``on the right'' of $\Hil$,
by $\rho(a)\bh=(ba)\hat{\;}$ $\forall b\in\MvN$.  Then $\alpha$
is left--right $\MvN$--equivariant if and only if
$(\pi(a)\otimes1)\alphah=(1\otimes\rho(a))\alphah$ $\forall a\in\MvN$.
So let us call $X\in\Hil\otimes\Hil$ left--right $\MvN$--equivariant
if $(\pi(a)\otimes1)X=(1\otimes\rho(a))X$ $\forall a\in\MvN$.
Then $\ef(\MvN,\phi)$ equals the infimum of
$\nm{(P_{\Hilo}\otimes P_{\Hilo})X}$ such that $X\in\Hil\otimes\Hil$
is left--right $\MvN$--equivariant and $\langle X,\xi\otimes\xi\rangle=1$.
We will use this to calculate expansion factors in~\S3.
\endproclaim

\proclaim{Lemma 2.3}
Let $\MvN$ be a diffuse von Neumann algebra, {\rm i.e\.} having no
nonzero minimal projections, and let $\phi$ be any n.f\.
state on $\MvN$.
Then every left--right $\MvN$--equivariant element of $\Hil\otimes\MvN$
equals zero, so $\ef(\MvN,\phi)=+\infty$.
\endproclaim
\demo{Proof}
Let $\alpha\in\Hil\otimes\MvN$ be left--right $\MvN$--equivariant,
and let $\Ac$ be a maximal abelian subalgebra (MASA) in $\MvN$.
It is well known that there is a conditional expectation of norm~$1$,
$E:\MvN\rightarrow\Ac$ (see~\cite{S}, 10.15).
Let $\HilA=L^2(\Ac,\phi\restrict_\Ac)\subseteq\Hil$
and $E\alpha=E\circ\alpha\restrict_{\HilAbar}:\HilAbar\rightarrow\Ac$.
Although $E\alpha$ is a bounded linear function and
$\nm{E\alpha}\le\nm{\alpha}$, in general $E\alpha\not\in\HilA\otimes\Ac$.
However $E\alpha$ is left--right $\Ac$--equivariant.
Indeed, we have for $a\in\Ac$ and $\eta\in\HilA$ that
$(E\alpha)(\pibar(a)\eta)=E(\alpha(\pibar(a)\eta))=E(\alpha(\eta)a)
=E(\alpha(\eta))a=(E\alpha(\eta))a$.
We will show that $E\alpha$ must be zero.
Let us take $\Ac=L^\infty([0,1])$, $\HilA=L^2([0,1])$ with
respect to Lebesque measure, and denote $E\alpha=T$.
Thus $T:\overline{L^2([0,1])}\rightarrow\L^\infty([0,1])$
is bounded linear.
The fact that $E\alpha$ is left--right $\Ac$--equivariant gives that
$$ T(fg)=\overline fT(g)\;\;\forall\,f\in L^\infty([0,1]),\,g\in L^2([0,1]). $$
So $T$ is determined by its value at the constant function $1$.
Suppose for contradiction that  $T\neq0$.
Then $T(1)\neq0$.
Let $\delta>0$ and $S\subseteq[0,1]$ be measurable of nonzero measure
such that $|T(1)(t)|>\delta$ $\forall t\in S$.
Without loss of generality we may suppose $S\supseteq[0,R]$, some $0<R\le1$.
For $n\ge1$ let $f_n\in\L^\infty([0,1])$ be
$$ f_n(t)=\cases t^{-1/4} & \text{ if }\frac1n\le t\le1 \\
                 0        & \text{ otherwise.}
\endcases $$
Then $\nm{f_n}_{L^2}^2=\int_{1/n}^1t^{-1/2}\dif t<2$.
But $\nm{T(f_n)}_{L^\infty}=\nm{f_nT(1)}_{L^\infty}\ge n^{1/4}\delta$
for $n\ge\frac1R$.
This contradicts the boundedness of $T$.
So $T$ must equal $0$.

  Hence we have shown that for every MASA $\Ac\subseteq\MvN$ and every
conditional expectation $E:\MvN\rightarrow\Ac$ we must have
$E\circ\alpha\restrict_{\HilAbar}=0$.
Suppose for contradiction that $\alpha\neq0$.
Then since $\pibar(\MvN)\xi$ is dense in $\Hilbar$ and since
$\alpha$ is left--right $\MvN$--equivariant, we must have
$\alpha(\xi)=x\neq0$.
Suppose the real part of $x$ is nonzero, $\xRe=\frac{x+x^*}2\neq0$,
let $\Ac$ be a MASA in $\MvN$ that contains
all the spectral projections of $\xRe$ and let
$E:\MvN\rightarrow\Ac$ be a conditional expectation.
Thus $E\alpha(\xi)=E(x)=\xRe+iE(\xIm)$.
Since $E$ is positive we know that $E(\xIm)$ is self--adjoint and hence that
$E(x)\neq0$.
This contradicts that $E\alpha=0$.
A similar argument works if $\xIm\neq0$.
\QED

\proclaim{Lemma 2.4}
Let $(\MvN,\phi)$ be a von Neumann algebra with n.f\. state
and suppose that $1\in\NvN\subseteq\MvN$ is a W$^*$--subalgebra
and $E:\MvN\rightarrow\NvN$ is a projection of norm $1$ satisfying
$\phi\circ E=\phi$.
Then $\ef(\MvN,\phi)\ge\ef(N,\phi\restrict_\NvN)$.
\endproclaim
\demo{Proof}
Let $\HilN=L^2(\NvN,\phi\restrict_\NvN)\subseteq\Hil$ and
$P_{\HilN}$ the orthogonal projection of $\Hil$ onto $\HilN$.
Since $E$ preserves $\phi$ it follows that
$(Ex)\hat{\;}=P_{\HilN}(\xh)$ $\forall x\in\MvN$.
Suppose $\alpha\in\Hil\otimes\MvN$ is left--right $\MvN$--equivariant.
Then letting $E\alpha=E\circ\alpha\restrict_{\HilNbar}:\HilNbar\rightarrow\NvN$
we have as in the proof of Lemma~2.3
that $E\alpha$ is left--right $\NvN$--equivariant.
We would like to show that for $E\alpha\in\HilN\otimes\NvN$.
First note that for $x\in\MvN$ we have
$E(x)\restrict_{\HilN}=(P_{\HilN}\circ x)\restrict_{\HilN}$.
Hence for $v\in\HilN$, $(E\alpha)(\cdot)v=(P_{\HilN}\circ\alpha(\cdot))v$, so
$$ ((E\alpha)(\cdot)v)\restrict_{\HilNbar}
=(P_{\HilN}\otimes P_{\HilN})
\bigl((\alpha(\cdot)v)\restrict_{\Hilbar_\NvN}\bigr)
\in\HilN\otimes\HilN. \tag6 $$
Thus $E\alpha\in\HilN\otimes\NvN$ and $\nm{E\alpha}_{HS}\le\nm{\alpha}_{HS}$.
A particular case of~(6) is
$(E\alpha)\hat{\;}=(P_{\HilN}\otimes P_{\HilN})\alphah$.
Since $E\alpha$ is left--right $\NvN$--equivariant we have
$\nm{(P_{\HilNo}\otimes P_{\HilNo})\alphah}\ge\ef(\NvN,\phi\restrict_\NvN)
|\langle(E\alpha)\hat{\;},\xi\otimes\xi\rangle|$,
where of course $P_{\HilNo}$ is the orthogonal projection of $\Hil$
onto $\HilNo=\HilN\ominus\Cpx\xi$.
But $\HilNo=\HilN\bigcap\Hilo$, so
$\nm{(P_{\Hilo}\otimes P_{\Hilo})\alphah}
\ge\nm{(P_{\HilNo}\otimes P_{\HilNo})\alphah}
\ge\ef(\NvN,\phi\restrict_\NvN)
|\langle(E\alpha)\hat{\;},\xi\otimes\xi\rangle|
=\ef(\NvN,\phi\restrict_\NvN)
|\langle\alphah,\xi\otimes\xi\rangle|$.
\QED

\proclaim{Theorem 2.5}
Let $(\MvN_\iota,\phi_\iota)$, ($\iota=1,2$) be von Neumann algebras
with faithful normal
states, and let $(\MvN,\phi)=(\MvN_1,\phi_1)*(\MvN_2,\phi_2)$.
Suppose that the linear dimensions satisfy
$\dim\MvN_\iota\ge2$, ($\iota=1,2$) and $\dim\MvN_1+\dim\MvN_2\ge5$,
and that $\ef(\MvN_1,\phi_1)\ef(\MvN_2,\phi_2)\ge1$,
(where $0\cdot(+\infty)=0$).
Then $\MvN$ is a factor and
$$ T(\MvN)=I(\MvN,\phi)=I(\MvN_1,\phi_1)\cap I(\MvN_2,\phi_2).
\tag7 $$
\endproclaim

The proof follows from the following proposition.

\proclaim{Proposition 2.6}
Under the hypotheses of the above theorem, if $u\in\MvN$ is a unitary such
that $\phi\circ\Ad_u=\phi$ ({\it i.e\.} u is in the centralizer, $\MvN_\phi$)
and $\Ad_u(\MvN_\iota)\subseteq\MvN_\iota$ ($\iota=1,2$), then $u\in\Cpx1$.
\endproclaim
\demo{Proof}
Let us use the notation of \S1 with index set $I=\{1,2\}$.
Let $AS=\{s=(s_1,\ldots,s_n)\in\bigcup_{m\ge1}\{1,2\}^m
\mid s_j\neq s_{j+1}\forall1\le j\le n-1\}$ be the set of all
alternating sequences of $1$'s and $2$'s.
For $s=(s_1\,\ldots,s_n)\in AS$ let
$\Hilo_s=\Hilo_{s_1}\otimes\cdots\otimes\Hilo_{s_n}\subseteq\Hil$
and let $P_{\Hilo_s}$ denote the orthogonal projection of $\Hil$
onto $\Hilo_s$.
Suppose that $u\not\in\Cpx1$ and consider $\uh\in\Hil$.
We will obtain a contradiction by showing that $\nm{P_{\Hilo_s}(\uh)}$
as $s$ increases in length are too large for $\uh$ to have bounded
norm in $\Hil$.

  Let $s=(s_1,\ldots,s_n)\in AS$ and suppose $s_1=s_n=k'$
and let $k\in\{1,2\}\backslash\{k'\}$.
Suppose $\gamma\in\Hilo_s$, $\nm{\gamma}=1$.
We define
$$ \Hil_k\otimes\gamma\otimes\Hil_k
=\Cpx\gamma\oplus(\Hilo_k\otimes\gamma)\oplus(\gamma\otimes\Hilo_k)
\oplus(\Hilo_k\otimes\gamma\otimes\Hilo_k), $$
where $\Hilo_k\otimes\gamma$ is the obvious subspace of
$\Hilo_k\otimes\Hilo_s=\Hilo_{(k,s)}$, {\it etcetera}.
$\Hil_k\otimes\gamma\otimes\Hil_k$ is isometrically identified
with $\Hil_k\otimes\Hil_k$ by identifying
$$ \alignat 3
& \gamma &\quad& \text{ with }\quad&\xi_k & \otimes\xi_k \\
\Hilo_k\otimes & \gamma &\quad& \text{ with }\quad&\Hilo_k & \otimes\xi_k \\
& \gamma\otimes\Hilo_k & \quad&\text{ with }\quad&\xi_k & \otimes\Hilo_k \\
\Hilo_k\otimes & \gamma\otimes\Hilo_k &\quad& \text{ with }\quad
&\Hilo_k & \otimes\Hilo_k.
\endalignat $$
Let $Q_\gamma:\Hil\rightarrow\Hil_k\otimes\Hil_k$ be the orthogonal projection
from $\Hil$ onto $\Hil_k\otimes\gamma\otimes\Hil_k$ followed by the
above identification of $\Hil_k\otimes\gamma\otimes\Hil_k$ with
$\Hil_k\otimes\Hil_k$.
One sees immediately that for $a\in\MvN_k$ and $\zeta\in\Hil$ one has
$Q_\gamma(\pi(a)\zeta)=(\pi(a)\otimes1)Q_\gamma(\zeta)$.

\proclaim{Lemma 2.7}
For $a\in\MvN$ and $\gamma$ as above there is a unique
$F_{\gamma,a}\in\Hil_k\otimes\MvN_k$ given by
$$ \langle F_{\gamma,a}(\eta)v,w\rangle_{\Hil_k}
=\langle av,Q_\gamma^*(\eta\otimes w)\rangle_{\Hil}
=\langle Q_\gamma(av),\eta\otimes w\rangle_{\Hil_k\otimes\Hil_k}\;\;
\forall\,\eta\in\Hilbar_k\text{ and }v,w\in\Hil_k, \tag8 $$
identifying $\Hilbar_k$ and $\Hil_k$ as sets and viewing $\Hil_k\subseteq\Hil$.
\endproclaim
\demo{Proof}
Clearly~(8) defines a bounded operator $F_{\gamma,a}(\eta)$ on
$\Hil_k$ for every $\eta\in\Hilbar_k$.
Also for
every fixed $v\in\Hil$, $F_{\gamma,a}(\cdot)v$ is Hilbert--Schmidt
and $\nm{F_{\gamma,a}(\cdot)v}_{HS}\le\nm{av}$.
So it suffices to show that $F_{\gamma,a}(\eta)\in\MvN_k$ for every
$\eta\in\Hilbar_k$ and every $a\in\MvN$.
Consider a word $b=b_1b_2\cdots b_m\in\MvNo_{t_1}\MvNo_{t_2}\cdots\MvNo_{t_m}$,
where $t=(t_1,\ldots,t_m)\in AS$.
Then $F_{\gamma,b}$ is a simple tensor in $\Hil_k\otimes\MvN_k$.
Indeed,
$$ F_{\gamma,b}=\cases
\langle \bh_1\otimes\cdots\otimes\bh_m,\gamma\rangle\quad\xi\otimes1
 &\text{ if }t=s \\
\langle \bh_2\otimes\cdots\otimes\bh_m,\gamma\rangle\quad\bh_1\otimes1
 &\text{ if }t=(k,s) \\
\langle \bh_1\otimes\cdots\otimes\bh_{m-1},\gamma\rangle\quad\xi\otimes b_m
 &\text{ if }t=(s,k) \\
\langle \bh_2\otimes\cdots\otimes\bh_{m-1},\gamma\rangle\quad\bh_1\otimes b_m
 &\text{ if }t=(k,s,k) \\
0 &\text{ otherwise.}
\endcases $$
Hence if $a$ is a linear combination of such words, then
$F_{\gamma,a}\in\MvN_k$.
However the $*$--subalgebra of such linear combinations of words is
dense in $\MvN$, so an arbitrary $a\in\MvN$ is the weak operator limit of a net
$(a_\lambda)_{\lambda\in\Lambda}$ of them,
and as is clear from~(8), $F_{\gamma,a}(\eta)$ is then the weak operator
limit of $(F_{\gamma,a_\lambda}(\eta))_{\lambda\in\Lambda}$, implying that
$F_{\gamma,a}\in\MvN_k$.
This proves Lemma~2.7.
\QED

Continuing with the proof of Proposition~2.6, let
$\tau_\iota=\Ad_u\restrict_{\MvN_\iota}$ ($\iota\in\{1,2\}$).
Since $\phi\restrict_{\MvN_\iota}\circ\tau_\iota=\phi\restrict_{\MvN_\iota}$,
there is a corresponding unitary $V_\iota\in\Bof(\Hil_\iota)$ given by
$V_\iota\bh=(\tau_\iota(b))\hat{\;}$ $\forall b\in\MvN_\iota$.
Note that $V_\iota\xi_\iota=\xi_\iota$.
One easily checks that $\pi_\iota(\tau_\iota(a))=V_\iota\pi_\iota(a)V_\iota^*$
$\forall a\in\MvN_\iota$.
\proclaim{Lemma 2.8}
Let $s$ and $\gamma$ be as above.
Let $G_\gamma=\tau_k\circ F_{\gamma,u}:\Hilbar_k\rightarrow\MvN_k$.
Then $G_\gamma\in\Hil_k\otimes\MvN_k$ and $G_\gamma$ is left--right
$\MvN_k$--equivariant.
\endproclaim
\demo{Proof}
Clearly $G_\gamma$ is bounded.
Note that $\nm{G_\gamma(\eta)v}=\nm{V_kF_{\gamma,u}(\eta)V_k^*v}$,
so $\nm{G_\gamma(\cdot)v}_{HS}=\nm{F_{\gamma,u}(\cdot)V_k^*v}_{HS}$
and hence $G_\gamma\in\Hil_k\otimes\MvN_k$ and
$\nm{G_\gamma}_{HS}=\nm{F_{\gamma,u}}_{HS}$.
To show that $G_\gamma$ is left--right $\MvN_k$--equivariant,
for $a\in\MvN_k$, $\eta\in\Hilbar_k$ and $v,w\in\Hil_k$ we have
$$ \align
\langle G_\gamma(\pibar(a)\eta)v,w\rangle
&=\langle F_{\gamma,u}(\pibar(a)\eta)V_k^*v,V_k^*w\rangle
=\langle u(V_k^*v),Q_\gamma^*(\pi(a^*)\eta\otimes(V_k^*w))\rangle
=\langle u(V_k^*v),\pi(a^*)Q_\gamma^*(\eta\otimes V_k^*w)\rangle= \\
&=\langle au(V_k^*v),Q_\gamma^*(\eta\otimes V_k^*w)\rangle
=\langle (u\tau^{-1}_k(a))(V_k^*v),Q_\gamma^*(\eta\otimes V_k^*w)\rangle
=\langle uV_k^*(av),Q_\gamma^*(\eta\otimes V_k^*w)\rangle=\\
&=\langle F_{\gamma,u}(\eta)V_k^*(av),V_k^*w\rangle
=\langle G_\gamma(\eta)av,w\rangle,
\endalign $$
so $G_\gamma(\pibar(a)\eta)=G_\gamma(\eta)a$,
{\it i.e\.} $G_\gamma$ is left--right $\MvN_k$--equivariant.
This proves Lemma~2.8.
\QED

\proclaim{Lemma 2.9}
There is $s=(s_1,\ldots,s_n)\in AS$ with $s_1=s_n$ such that
$P_{\Hilo_s}(\uh)\neq0$.
\endproclaim
\demo{Proof}
Since we supposed that $u\not\in\Cpx1$, there is $t=(t_1,\ldots,t_n)\in AS$
such that $P_{\Hilo_t}(\uh)\neq0$.
If $t_1=t_n$ we are done.
If $t_1\neq t_n$ let $k=t_1$, $k'=t_n$.
We may suppose $\dim\MvN_k\ge3$.
(If instead $\dim\MvN_{k'}\ge3$ the same argument but reflected
will work.)
Then there is $\gamma\in\Hilo_{t'}$, where $t'=(t_2,\ldots,t_n)$,
such that
$$ Q_\gamma(\uh)=r(\xi_k\otimes\xi_k)+\zeta_1\otimes\xi_k
+\xi_k\otimes\zeta_2+X, $$
where $r\in\Cpx$, $\zeta_1,\zeta_2\in\Hilo_k$, $X\in\Hilo_k\otimes\Hilo_k$
and $\zeta_1\neq0$.
To prove the lemma it will suffice to show that $X\neq0$, because then
$P_{\Hilo_{(t,k)}}(\uh)\neq0$.
Suppose for contradiction that $X=0$.
Then
$$ \langle G_\gamma(\eta)\xi_k,w\rangle
=\langle Q_\gamma(\uh),\eta\otimes V_k^*w\rangle
=r\langle\xi_k,\eta\rangle\langle\xi_k,w\rangle
+\langle\zeta_1,\eta\rangle\langle\xi_k,w\rangle
+\langle\xi_k,\eta\rangle\langle V_k\zeta_2,w\rangle, $$
so $G_\gamma(\eta)\xi_k=\langle r\xi_k+\zeta_1,\eta\rangle\xi_k
+\langle\xi_k,\eta\rangle V_k\zeta_2\in(\MvN_k)\hat{\;}$
$\forall\eta\in\Hilbar_k$.
Hence $V_k\zeta_2=\bh$, some $b\in\MvN_k$ and
$$ G_\gamma=(r\xi_k+\zeta_1)\otimes1+\xi_k\otimes b. $$
{}From the left--right $\MvN_k$--equivariance of $G_\gamma$ we have that
$$ a(r\xi_k+\zeta_1)\otimes1+\ah\otimes b
=(r\xi_k+\zeta_1)\otimes a+\xi_k\otimes ba\;\;\forall a\in\MvN_k. \tag9 $$
If $b=0$ then we can get a contradiction to~(9) by choosing any
$a\not\in\Cpx1$, (remembering that
$\zeta_1\neq0$ and $\zeta_1\perp\xi_k$).
If $b\neq0$ then since $\dim\MvN_k\ge3$ there is $a\in\MvN_k$
such that $\ah\perp\xi_k$ and $\ah\perp\zeta_1$.
For this $a$, we have that the right hand side of~(9) is an element of
$\Hil_k\otimes\MvN_k$ sending $\ah\in\Hilbar_k$ to zero while
the left hand side sends $\ah$ to some multiple of $1$ plus
$\langle\ah,\ah\rangle b$, which is nonzero, a contradiction.
Hence $X$ must be nonzero and Lemma~2.9 is proved.
\QED

\proclaim{Lemma 2.10}
Let $s=(s_1,\ldots,s_n)\in AS$ with $s_1=s_n=k'$.
Let $k\in\{1,2\}\backslash\{k'\}$ and $t=(k,s_1,s_2,\ldots,s_n,k)$.
Then $\nm{P_{\Hilo_t}(\uh)}\ge\ef(\MvN_k,\phi_k)\nm{P_{\Hilo_s}(\uh)}$.
\endproclaim
\demo{Proof}
Let $\{\gamma_j\}_{j\in J}$ be an orthonormal basis for $\Hilo_s$.
Then
$$ \aligned
\nm{P_{\Hilo_s}(\uh)}^2
 &=\sum_{j\in J}|\langle\gamma_j,\uh\rangle|^2
 =\sum_{j\in J}|\langle\xi_k\otimes\xi_k,Q_\gamma(\uh)\rangle|^2 \\
\nm{P_{\Hilo_t}(\uh)}^2
 &=\sum_{j\in J}\nm{(P_{\Hilo_k}\otimes P_{\Hilo_k})Q_\gamma(\uh)}^2.
\endaligned \tag10 $$
But note that for $\gamma\in\{\gamma_j\}_{j\in J}$,
$$
\langle Q_\gamma(\uh),\xi_k\otimes\xi_k\rangle
=\langle u(V_k^*\xi_k),Q_\gamma^*(\xi_k\otimes V_k^*\xi_k\rangle
=\langle G_\gamma(\xi_k)\xi_k,\xi_k\rangle
=\langle\Gh_\gamma,\xi_k\otimes\xi_k\rangle,
$$
so $\nm{P_{\Hilo_s}(\uh)}^2
=\sum_{j\in J}|\langle\Gh_\gamma,\xi_k\otimes\xi_k\rangle|$.
Let $\{v_\lambda\}_{\lambda\in\Lambda}$ be an orthonormal basis for $\Hilo_k$.
Then for $\gamma\in\{\gamma_j\}_{j\in J}$,
$$ \aligned \dsize
\nm{(P_{\Hilo_k}\otimes P_{\Hilo_k})Q_\gamma(\uh)}^2
&=\sum_{\lambda,\mu\in\Lambda}
 |\langle u(V_k^*\xi_k),Q_\gamma^*(v_\lambda\otimes V_k^*v_\mu)\rangle|^2
=\sum_{\lambda,\mu\in\Lambda}
 |\langle G_\gamma(v_\lambda)\xi_k,v_\mu\rangle|^2 \\ \dsize
&=\sum_{\lambda,\mu\in\Lambda}
 |\langle\Gh_\gamma,v_\lambda\otimes v_\mu\rangle|^2
=\nm{(P_{\Hilo_k}\otimes P_{\Hilo_k})(\Gh_\gamma)}^2.
\endaligned $$
But since $G_\gamma$ is left--right $\MvN_k$--equivariant, we have by
definition that
$$ \nm{(P_{\Hilo_k}\otimes P_{\Hilo_k})(\Gh_\gamma)}^2
\ge\ef(\MvN_k,\phi_k)^2|\langle \Gh_\gamma,\xi_k\otimes\xi_k\rangle|^2. $$
Now use this together with~(10) to prove Lemma~2.10.
\QED

\demo{Proof of Proposition 2.6}
Lemma~2.9 guarantees us that there is $s=(s_1,\ldots,s_n)\in AS$ such
that $s_1=s_n=k'$ and $P_{\Hilo_s}(\uh)\neq0$.
Let $k\in\{1,2\}\backslash\{k'\}$ and
$$ t(m)=(\underbrace{k',k,\ldots,k',k}_{m\text{ times }k',k},
s_1,\ldots,s_n,\underbrace{k,k',\ldots,k,k'}_{m\text{ times }k,k'}). $$
Then by Lemma~2.10,
$$ \nm{P_{\Hilo_{t(m)}}(\uh)}^2
\le\ef(\MvN_{k'},\phi_{k'})^{2m}\ef(\MvN_k,\phi_k)^{2m}
\nm{P_{\Hilo_s}(\uh)}^2. $$
Under the hypotheses of Theorem~2.5,
this gives a contradiction to $\uh\in\Hil$.
\QED

\demo{Proof of Theorem 2.5}
To prove the factoriality of $\MvN$ it suffices to show that any unitary $u$
that is in the center of $\MvN$ must be a constant.
However, such a unitary satisfies the hypotheses of Proposition~2.6,
so is constant.
To prove~(7), suppose $t\in T(\MvN)$, {\it i.e\.}
$\sigma_t=\Ad_u$, $u\in\MvN$.
Of course $\sigma_t$ preserves $\phi$, and by Theorem~1 we know that
$\sigma_t$ sends $\MvN_\iota$ into $\MvN_\iota$ ($\iota=1,2$), so again
by Proposition~2.6 we must have that $u$ is constant, so that
$t\in I(\MvN,\phi)$.
Thus $T(\MvN)=I(\MvN,\phi)$.
But also by Theorem~1 we have the last equality in~(7).
\QED

\noindent{\bf\S3 Calculations.}

  In this section we calculate the expansion factors of various algebras.
Actually, the methods used here allow one, with a bit of work, to calculate
$\ef(\MvN,\phi)$ for any pair $(\MvN,\phi)$.
But since Theorem~2.5 is sufficient but not necessary for factoriality,
(see Remark~4.1), It would not be rewarding to compute $\ef(\MvN,\phi)$
for the most general $(\MvN,\phi)$.
We compute enough special cases to show (Lemma~3.2) that the hypotheses
of Theorem~2.5 are satisfied in gratifyingly many cases, and also
to illustrate some of the limitations of Theorem~2.5, (see Remark~4.1).

\proclaim{Lemma 3.1}
Let $A$ be a finite dimensional commutative algebra with faithful state
$\phi$.
Suppose $e_1,\ldots,e_n$ are the minimal projections of $A$ having
traces $\phi(e_\iota)=\lambda_\iota$, ($1\le\iota\le n$), with
$\lambda_1\ge\lambda_2\ge\cdots\ge\lambda_n>0$.
Then
$$ \ef(A,\phi)^2=1+\biggl(\sum_{\iota=1}^n\frac{\lambda_\iota^2}
{1-2\lambda_\iota}\biggr)^{-1}, \tag11 $$
where if $\lambda_1=\frac12$ we take the right hand side of~(11) to be $1$.
\endproclaim
\demo{Proof}
We will use the method described in Remark~2.2 to calculate the expansion
factor.
$\Hil=L^2(A,\phi)$ has orthogonal basis $\{\eh_\iota\mid1\le\iota\le n\}$.
Hence every $X\in\Hil\otimes\Hil$ can be written $X=\sum_{\iota,j=1}^n
c_{\iota,j}\eh_\iota\otimes \eh_j$.
Since for $X$ to be left--right $A$--invariant
means $(\pi(e_p)\otimes1)X=(1\otimes\rho(e_p))X$, we see that
$c_{\iota,j}=0$ if $\iota\neq j$ and that
$X$ is left--right
$A$--invariant if and only if it is of the form
$$ X=\sum_{\iota=1}^n c_\iota \eh_\iota\otimes \eh_\iota, $$
for $c_\iota\in\Cpx$.
Now in $\Hil$, $\nm{\eh_\iota}^2=\phi(e_\iota e_\iota)=\lambda_\iota$
and $P_\xi(\eh_\iota)=\lambda_\iota\xi$, so
$$ \alignat 2
(P_\xi\otimes P_\xi)(X)
&=\biggl(\sum_{\iota=1}^n c_\iota\lambda_\iota^2\biggr)\xi\otimes\xi,
&\qquad\nm{(P_\xi\otimes P_\xi)(X)}^2
&=\biggl|\sum_{\iota=1}^n c_\iota\lambda_\iota^2\biggr|^2, \\
(P_\xi\otimes I)(X)
&=\sum_{\iota=1}^n c_\iota\lambda_\iota(\xi\otimes\eh_\iota),
&\nm{(P_\xi\otimes I)(X)}^2
&=\sum_{\iota=1}^n |c_\iota|^2\lambda_\iota^3, \\
(I\otimes P_\xi)(X)
&=\sum_{\iota=1}^n c_\iota\lambda_\iota(\eh_\iota\otimes\xi),
&\nm{(I\otimes P_\xi)(X)}^2
&=\sum_{\iota=1}^n |c_\iota|^2\lambda_\iota^3, \\
&&\nm{X}^2
&=\sum_{\iota=1}^n |c_\iota|^2\lambda_\iota^2.
\endalignat $$
Thus
$$ \align
 \nm{(P_\Hilo\otimes P_\Hilo)(X)}^2&=\nm X^2
-\nm{(P_\xi\otimes I)(X)}^2
-\nm{(I\otimes P_\xi)(X)}^2
+\nm{(P_\xi\otimes P_\xi)(X)}^2 \\
&=\sum_{\iota=1}^n |c_\iota|^2(\lambda_\iota^2-2\lambda_\iota^3)
+\biggl|\sum_{\iota=1}^n c_\iota\lambda_\iota^2\biggr|^2. \tag12
\endalign $$
As remarked in~2.2, $\ef(A,\phi)^2$ is the infimum of~(12)
as $(c_\iota)_{1\le\iota\le n}$ varies over $\Cpx^n$,
subject to
the constraint
$$ C=\sum_{\iota=1}^n c_\iota\lambda_\iota^2=1, \tag13 $$
{\it i.e.} the infimum
of
$$ V=1+\sum_{\iota=1}^n |c_\iota|^2\lambda_\iota^2(1-2\lambda_\iota)
\tag14 $$
subject to~(13).

  When $\lambda_1\le\frac12$, clearly $V\ge1$ everywhere,
and for the special case $\lambda_1=\frac12$ we can chose $c_1=4$ and
$c_\iota=0$ for $2\le\iota\le n$ to give $V=1$, so the infimum
equals $1$, as required.  We may thus assume henceforth that
$\lambda_1\neq\frac12$.

  Let us show that it suffices to find the infimum of $V$ subject
to~(13) for $(c_\iota)_{1\le\iota\le n}\in\Real^n$.
Letting $c_\iota=a_\iota+ib_\iota$ for $a_\iota,b_\iota\in\Real$,
it suffices to show that
$V'=\sum_{\iota=1}^n b_\iota^2\lambda_\iota^2(1-2\lambda_\iota)$
is always non-negative when
$C'=\sum_{\iota=1}^n b_\iota\lambda_\iota^2=0$.
But if some such choice of $(b_\iota)_{1\le\iota\le n}$ gives
negative $V'$, then multiples of them also satisfy $C'=0$
and give $V'$ as large and negative as we please,
enough to make $V$ itself negative.
This is a contradiction, since $V$, being the square--norm of
a vector, must be nonnegative.
Hence we may assume $c_\iota\in\Real$ $\forall\iota$.

  In order to use the method of Lagrange multipliers, we need to know that
the infimum of $V$ subject to~(13) occurs at a relative minimum of $V$
on the manifold defined by~(13).  For $\lambda_1<\frac12$, this is obvious,
since it is then clear that $V\rightarrow+\infty$ as
$(c_\iota)_{1\le\iota\le n}\rightarrow\infty$.
For $\lambda_1>\frac12$ we proceed as follows.
For $s\in\Real$ let $V(s)$ denote the minimum value of $V$ subject to~(13)
with the value of $c_1$ fixed to be $s$.
We will show that $V(s)\rightarrow+\infty$ as $s\rightarrow\pm\infty$.
It is clear that with $c_1$ fixed, $V\rightarrow+\infty$ as
$(c_\iota)_{2\le\iota\le n}\rightarrow\infty$, so we may use the
method of Lagrange multipliers to find $V(s)$.
We have
$$ \frac{\partial C}{\partial c_\iota}=\lambda_\iota^2 \tag15 $$
and
$$ \frac{\partial V}{\partial c_\iota}
=2c_\iota\lambda_\iota^2(1-2\lambda_\iota), \tag16 $$
so the value of $V(s)$ occurs
where
$$ c_\iota=\frac r{2(1-2\lambda_\iota)},\;(2\le\iota\le n), $$
for some $ r\in\Real$.
{}From~(13) we have that
$$  r=(1-s\lambda_1^2)\biggl(\sum_{\iota=2}^n
\frac{\lambda_\iota^2}{2(1-2\lambda_\iota)}\biggr)^{-1}, $$
so
$$ \align
V(s)&=1+s^2\lambda_1^2(1-2\lambda_1)+ r^2
\sum_{\iota=2}^n \frac{\lambda_\iota^2(1-2\lambda_\iota)}
{4(1-2\lambda_\iota)^2} \\
&=1+s^2\lambda_1^2(1-2\lambda_1)+\tfrac12(1-s\lambda_1^2)^2
\biggl(\sum_{\iota=2}^n \frac{\lambda_\iota^2}
{2(1-2\lambda_\iota)}\biggr)^{-1} \\
&=s^2\lambda_1^2(1-2\lambda_1)+\tfrac12s^2\lambda_1^4
\biggl(\sum_{\iota=2}^n \frac{\lambda_\iota^2}
{2(1-2\lambda_\iota)}\biggr)^{-1}+\text{ lower terms in $s$.}
\endalign $$
We must show that
$$ \tfrac12\lambda_1^4
\biggl(\sum_{\iota=2}^n \frac{\lambda_\iota^2}
{2(1-2\lambda_\iota)}\biggr)^{-1}
>\lambda_1^2(2\lambda_1-1), $$
{\it i.e\.}
$$ \frac{\lambda_1^2}{2\lambda_1-1}
>\sum_{\iota=2}^n \frac{\lambda_\iota^2}{1-2\lambda_\iota}. \tag17 $$
Consider the function $f(t)=\frac{t^2}{1-2t}$.
One easily verifies that $f(0)=0$ and that on $(0,\frac12)$,
$f'(t)$ and $f''(t)$ are strictly positive.
Thus if $t_1,t_2\in(0,\frac12)$, $t_1\le t_2$ and $t_1+t_2<\frac12$,
then the point $(t_1+t_2,f(t_1+t_2))$ lies above the line
$y=\frac{f(t_2)}{t_2}x$,
and the point $(t_1,f(t_1))$ lies below this
line.
Hence
$$ f(t_1)+f(t_2)\le\frac{f(t_2)}{t_2}t_1+f(t_2)\le f(t_1+t_2). $$
Consequently we have
$$ \sum_{\iota=2}^n\frac{\lambda_\iota^2}{1-2\lambda_\iota}
=\sum_{\iota=2}^n f(\lambda_\iota)
\le f\biggl(\sum_{\iota=2}^n\lambda_\iota\biggr)
=f(1-\lambda_1)
=\frac{(1-\lambda_1)^2}{2\lambda_1-1}. $$
Thus to show~(17) it suffices to show
$$\frac{\lambda_1^2}{2\lambda_1-1}
>\frac{(1-\lambda_1)^2}{2\lambda_1-1}, $$
but since $\lambda_1>\frac12$, this is true.
We have thus shown that $V(s)$ is a quadratic polynomial in $s$ with
positive leading coefficient, so $V(s)\rightarrow+\infty$ as
$s\rightarrow\pm\infty$.

  We are now justified in using the method of Lagrange multipliers
to find the minimum of~(14) subject to the constraint~(13).
{}From~(15) and~(16) we see that the minimum value of $V$ occurs when
$$ c_\iota=\frac r{2(1-2\lambda_\iota)},\;(1\le\iota\le n) $$
for some $ r\in\Real$.
{}From~(13) we get
$$  r=\biggl(\sum_{\iota=1}^n\frac{\lambda_\iota^2}
{2(1-2\lambda_\iota)}\biggr)^{-1}, $$
so the value of $V$ at the global minimum is
$$ V_{\text{min}}=1+ r^2\sum_{\iota=1}^n
\frac{\lambda_\iota^2(1-2\lambda_\iota)}
{4(1-2\lambda_\iota)^2}
=1+\frac r2
=1+\biggl(\sum_{\iota=1}^n\frac{\lambda_\iota^2}{1-2\lambda_\iota}\biggr)
^{-1}. $$
\QED

\proclaim{Lemma 3.2}
Let $\MvN$ be a von Neumann algebra and $\phi$ a faithful state on $\MvN$.
If there are no minimal projections $e\in\MvN$ such that $\phi(e)>\frac12$,
then $\ef(\MvN,\phi)\ge1$.
\endproclaim
\demo{Proof}
We can write $\MvN=\Ac_1\oplus\Ac_2$, where $\Ac_2$ is diffuse and where
$\Ac_1$ has the property that every central projection contains a nontrivial
minimal projection.  (Thus $\Ac_1$ is a possibly infinite direct sum of
type~I factors.)  Let $\Hil_\iota=L^2(\Ac_\iota,\phi\restrict_{\Ac_\iota})$
($\iota=1,2$), $\Hil=L^2(\MvN,\phi)=\Hil_1\oplus\Hil_2$ and consider
$\alpha\in\Hil\otimes\MvN$ that is left--right $\MvN$--equivariant.
Let $p\in\MvN$ be the central projection $1\oplus0$.
Then since $\alpha(\pibar(p)\eta)=\alpha(\eta)p$
$\forall\eta\in\Hilbar=\Hilbar_1\oplus\Hilbar_2$, we have that
$\alpha=\alpha_1\oplus\alpha_2$, where
$\alpha_\iota:\Hil_\iota\rightarrow\Ac_\iota$ is left--right
$\Ac_\iota$--equivariant.
But $\alpha_2=0$ by Lemma~2.3.
Thus $(\MvN,\phi)$ has the same expansion factor as what we get
when we replace $\Ac_2$ by a diffuse commutative von Neumann algebra,
so we assume without loss of generality that $\Ac_2$ is commutative.
Then there is a commutative subalgebra $\Bc\subseteq\MvN$
having no minimal projections $e\in\Bc$ with $\phi(e)>\frac12$
and a conditional
expectation $E:\MvN\rightarrow\Bc$ such that $\phi\circ E=\phi$.
Indeed, since $\Ac_1$ is a direct sum of type~I factors, we need only
see that a type~I factor $\NvN$ with faithful state $\psi$
has a commutative subalgebra $\Cc\in\NvN$ and a $\psi$--preserving
conditional expectation of $\NvN$ onto $\Cc$.
However, $\psi(\cdot)=\Tr(h\cdot)$, where $\Tr$ is a faithful trace on $\NvN$
and where $h\in\NvN$ is a positive trace--class operator, each of
whose spectral projections is a finite projection.  Let $\Cc$ be the MASA
containing all the spectral projections of $h$.  Let $(e_\iota)_{\iota\in I}$
be the set of minimal projections of $\Cc$ and let $E:\NvN\rightarrow\Cc$
be $E(x)=\sum_{\iota\in I}\psi(e_\iota xe_\iota)e_\iota$.
Then $E$ is a $\psi$--preserving conditional expectation.
Thus we have shown the existence of commutative $\Bc$ as described above
and a $\phi$--preserving
conditional expectation.
So by Lemma~2.4, $\ef(\MvN,\phi)\ge\ef(\Bc,\phi\restrict_\Bc)$.
But $\Bc$ has a finite dimensional subalgebra $\Bc_0$
having no minimal projections $e$ with $\phi(e)>\frac12$, and
$\Bc_0$ is of course the image of a $\phi$--preserving conditional expectation
$E_0:\Bc\rightarrow\Bc_0$, so applying Lemma~2.4 again yields
$\ef(\MvN,\phi)\ge\ef(\Bc_0,\phi\restrict_{\Bc_0})$.
An application of Lemma~3.1 now completes the proof.
\QED

\proclaim{Lemma 3.3}
Let $A=M_n(\Cpx)$ be the algebra of $n\times n$ complex matrices,
containing a system of matrix units $\{e_{\iota j}\mid1\le\iota,j\le n\}$ and
equipped with a faithful state $\phi$ such that
$$ \phi(e_{\iota j})=\cases \lambda_\iota&\text{ if }\iota=j \\
0&\text{ if }\iota\neq j,
\endcases $$
where $\lambda_1\ge\lambda_2\ge\cdots\ge\lambda_n>0$.
Then
$$ \ef(A,\phi)^2=1+\biggl(\sum_{p=1}^n\frac{\lambda_p^3}
{1-2\lambda_p^2}\biggr)^{-1}, $$
where we interpret the right hand side of the above equation to be $1$
if $\lambda_1=\frac1{\sqrt{2}}$.
\endproclaim
\demo{Proof}
The proof goes pretty much exactly as that of Lemma~3.1.
The set $\{\eh_{\iota j}\mid1\le\iota,j\le n\}$ is an orthogonal basis
for $\Hil=L^2(A,\phi)$, and an arbitrary element $X\in\Hil\otimes\Hil$
is $X=\sum_{\iota,j,k,l}c_{\iota j k l}\eh_{\iota j}\otimes\eh_{kl}$.
For $X$ to be left--right $A$--equivariant means
$(\pi(e_{pq})\otimes1)X=(1\otimes\rho(e_{pq}))X$, {\it i.e\.}
$$\sum_{1\le j,k,l\le n}c_{q,j,k,l}\eh_{pj}\otimes\eh_{kl}
=\sum_{1\le\iota,j,k\le n}c_{\iota,j,k,p}\eh_{\iota j}\otimes\eh_{kq}. $$
Thus $c_{\iota,j,k,l}=0$ unless $\iota=l$ and $c_{p,j,k,p}=c_{q,j,k,q}$
$\forall1\le p,q\le n$, so $X\in\Hil\otimes\Hil$
is left--right equivariant if and only if it is of the form
$$ X=\sum_{j,k=1}^nc_{j,k}\biggl(\sum_{p=1}^n\eh_{pj}\otimes\eh_{kp}\biggr) $$
for $c_{j,k}\in\Cpx$.
In $\Hil$ we have
$\nm{\eh_{st}}^2=\phi(e_{st}^*e_{st})=\phi(e_{tt})=\lambda_t$
and $P_\xi(\eh_{pq})=\delta_{pq}\lambda_p\xi$,
so we get
$$ \alignat 2
(P_\xi\otimes P_\xi)(X)
&=\biggl(\sum_{p=1}^n c_{p,p}\lambda_p^2\biggr)\xi\otimes\xi,
&\qquad\nm{(P_\xi\otimes P_\xi)(X)}^2
&=\biggl|\sum_{p=1}^n c_{p,p}\lambda_p^2\biggr|^2, \\
(P_\xi\otimes I)(X)
&=\sum_{j,k=1}^n c_{j,k}\lambda_j(\xi\otimes \eh_{kj}),
&\nm{(P_\xi\otimes I)(X)}^2
&=\sum_{j,k=1}^n |c_{j,k}|^2\lambda_j^3, \\
(I\otimes P_\xi)(X)
&=\sum_{j,k=1}^n c_{j,k}\lambda_k(\eh_{kj}\otimes\xi),
&\nm{(I\otimes P_\xi)(X)}^2
&=\sum_{j,k=1}^n |c_{j,k}|^2\lambda_j\lambda_k^2,
\endalignat $$
$$ \nm{X}^2
=\sum_{j,k=1}^n |c_{j,k}|^2\sum_{p=1}^n\lambda_j\lambda_p
=\sum_{j,k=1}^n |c_{j,k}|^2\lambda_j,
$$
$$ \align
 \nm{(P_\Hilo\otimes P_\Hilo)(X)}^2&=\nm X^2
-\nm{(P_\xi\otimes I)(X)}^2
-\nm{(I\otimes P_\xi)(X)}^2
+\nm{(P_\xi\otimes P_\xi)(X)}^2 \\
&=\sum_{j,k=1}^n|c_{j,k}|^2\lambda_j(1-\lambda_j^2-\lambda_k^2)
+\biggl|\sum_{p=1}^nc_{p,p}\lambda_p^2\biggr|^2.
\endalign $$

  Now $\ef(A,\phi)^2$ will equal the infimum of
$V=1+\sum_{j,k=1}^n|c_{j,k}|^2\lambda_j(1-\lambda_j^2-\lambda_k^2)$
as $(c_{j,k})_{1\le j,k\le n}$ ranges over $\Cpx^{n^2}$
subject to the constraint
$C=\sum_{p=1}^nc_{p,p}\lambda_p^2=1$.
Note that if $j\neq k$ then since $\lambda_j+\lambda_k\le1$ we
have $\lambda_j^2+\lambda_k^2<1$.
Hence we may assume without loss of generality that $c_{j,k}=0$
if $j\neq k$,
and rewriting $c_{p,p}=c_p$, we want to find the infimum of
$$ V=1+\sum_{p=1}^n|c_p|^2\lambda_p(1-2\lambda_p^2) \tag18 $$
 as
$(c_p)_{1\le p\le n}$ ranges over $\Cpx^n$ subject to
the constraint
$$ C=\sum_{p=1}^nc_p\lambda_p^2=1. \tag19 $$

  When $\lambda_1\le\frac1{\sqrt2}$, clearly $V\ge1$ everywhere.
For the special case $\lambda_1=\frac1{\sqrt2}$ we get $V=1$ by setting
$c_1=2$ and $c_p=0$ for $2\le p\le n$.
Thus the infimum is $1$, as required, and we may henceforth
assume that $\lambda_1\neq\frac1{\sqrt2}$.

  We can easily show that we may restrict ourselves to considering
$(c_p)_{1\le p\le n}\in\Real^n$, exactly as we did in the proof of Lemma~3.1.

  Also as in the proof of Lemma~3.1, the infimum of $V$ on the manifold
defined by~(19) clearly occurs at a local minimum if
$\lambda_1<\frac1{\sqrt2}$,
but we need to argue further to show this if $\lambda_1>\frac1{\sqrt2}$.
Again let $V(s)$ for $s\in\Real$ be the minimum of $V$ on the manifold
defined by~(19) when $c_1$ is fixed to have the value $s$.
We may clearly use the method of Lagrange multipliers to find $V(s)$.
We have
$$ \frac{\partial C}{\partial c_p}=\lambda_p^2 \tag20 $$
and
$$ \frac{\partial V}{\partial c_p}=2c_p\lambda_p(1-2\lambda_p^2), \tag21 $$
so the value of $V(s)$ occurs where
$$ c_p=\frac{ r\lambda_p}{2(1-2\lambda_p^2)},\;\;(2\le p\le n), $$
for some $ r\in\Real$.
{}From~(19) we obtain
$$  r=(1-s\lambda_1^2)\biggl(\sum_{p=2}^n
\frac{\lambda_p^3}{2(1-2\lambda_p^2)}\biggr)^{-1}. $$
Thus
$$ \align
V(s)&=1+s^2\lambda_1(1-2\lambda_1^2)
+ r^2\sum_{p=2}^n\frac{\lambda_p^3(1-2\lambda_p^2)}
{4(1-2\lambda_p^2)^2} \\
&=1+s^2\lambda_1(1-2\lambda_1^2)
+\tfrac12 r^2\sum_{p=2}^n\frac{\lambda_p^3}
{2(1-2\lambda_p^2)} \\
&=1+s^2\lambda_1(1-2\lambda_1^2)
+\tfrac12(1-s\lambda_1^2)^2\biggl(\sum_{p=2}^n\frac{\lambda_p^3}
{2(1-2\lambda_p^2)}\biggr)^{-1}. \tag22
\endalign $$
Now~(22) is a quadratic polynomial in $s$, and the coefficient of
$s^2$ is
$$ \lambda_1(1-2\lambda_1^2)+\frac{\lambda_1^4}2
\biggl(\sum_{p=2}^n\frac{\lambda_p^3}
{2(1-2\lambda_p^2)}\biggr)^{-1}, $$
which we want to show is positive, {\it i.e\.}
$$ \frac{\lambda_1^3}2>(2\lambda_1^2-1)
\biggl(\sum_{p=2}^n\frac{\lambda_p^3}
{2(1-2\lambda_p^2)}\biggr), $$
{\it i.e\.}
$$ \frac{\lambda_1^3}{2\lambda_1^2-1}>\sum_{p=2}^n
\frac{\lambda_p^3}{1-2\lambda_p^2}. $$
Consider the functions $g(t)=\frac{t^3}{1-2t^2}$.
Just as we did in the proof of Lemma~3.1 for $f$, we
can easily show that $g(t_1)+g(t_2)\le g(t_1+t_2)$ when
$t_1,t_2$ and $t_1+t_2$ are in $(0,\frac12)$.
But then we have
$$ \sum_{p=2}^n\frac{\lambda_p^3}{1-2\lambda_p^2}
=\sum_{p=2}^ng(\lambda_p)
\le g\biggl(\sum_{p=2}^n\lambda_p\biggr)
=g(1-\lambda_1)
=\frac{(1-\lambda_1)^3}{1-2(1-\lambda_1)^2}. $$
So we would like to show that
$$ \frac{\lambda_1^3}{2\lambda_1^2-1}
>\frac{(1-\lambda_1)^3}{1-2(1-\lambda_1)^2}
\text{ for }\tfrac1{\sqrt2}<\lambda_1<1, $$
which works out to showing that the polynomial
$p(\lambda_1)=-2\lambda_1^4+4\lambda_1^3+\lambda_1^2
-3\lambda_1+1>0
\text{ for }\tfrac1{\sqrt2}<\lambda_1<1$.
One finds that the roots of $p(t)=0$ are
$t=\frac{1\pm\sqrt{4\pm\sqrt{17}}}2$, and thus the
real roots are approximately $-0.925$ and $1.925$.
Since $p(t)\rightarrow-\infty$ as $t\rightarrow\pm\infty$,
this shows that $p(\lambda_1)>0$ for $\tfrac1{\sqrt2}<\lambda_1<1$.
Thus we are justified in using the method of Lagrange multipliers
to find $V$.

  From~(20) and~(21) we have that
$$ c_p=\frac{ r\lambda_p}{2(1-2\lambda_p^2)}\quad(1\le p\le n) $$
for some $r\in\Real$,
and hence from~(19) we get
$$  r=\biggl(\sum_{p=1}^n\frac{\lambda_p^3}
{2(1-2\lambda_p^2)}\biggr)^{-1}. $$
Thus from~(18), the value of $V$ at its global minimum is
$$ V_{\text{min}}
=1+ r^2\sum_{p=1}^n\frac{\lambda_p^3}{4(1-2\lambda_p^2)}
=1+\frac r2
=1+\biggl(\sum_{p=1}^n\frac{\lambda_p^3}
{1-2\lambda_p^2}\biggr)^{-1}. $$
\QED

\proclaim{Examples 3.4}\rm
\roster
\item"i)" If $A=\Cpx^n$ ($n\ge2$) and $\phi$ is the state giving
each minimal projection weight $\frac1n$, then
$$ \ef(A,\phi)^2=1+\left(n\frac{\frac1{n^2}}{1-\frac2n}\right)^{-1}
=n-1. $$
\item"ii)" If $A=M_n(\Cpx)$ with $\tau$ the normalized trace, then
$$ \ef(A,\tau)^2=1+\left(n\frac{\frac1{n^3}}{1-\frac2{n^2}}\right)^{-1}
=n^2-1. $$
\item"iii)" If $A=\Cpx\bigoplus\Cpx$ and $\phi$ is the state assigning
to one of the minimal projections the value $\lambda>\frac12$, then
$$ \ef(A,\phi)^2=1+\left(\frac{\lambda^2}{1-2\lambda}
+\frac{(1-\lambda)^2}{2\lambda-1}\right)^{-1}=0.
$$
\item"iv)" If $A=M_2(\Cpx)$, if $\phi$
is as described in Lemma~3.3
and if we write $\lambda=\lambda_1\in[\frac12,1)$, then
$$\ef(A,\phi)^2=1+\left(\frac{\lambda^3}{1-2\lambda^2}+\frac{(1-\lambda)^3}
{1-2(1-\lambda)^2}\right)^{-1}
=\frac{\lambda+\lambda^2-4\lambda^3+2\lambda^4}
{1-3\lambda+\lambda^2+4\lambda^3-2\lambda^4}.
$$
\endroster
\endproclaim

\noindent{\bf \S4. Consequences and limitations of the main theorem.}

\nopagebreak

\proclaim{Remark 4.1} Factoriality. \rm
\roster
\item"i)" From the calculations in \S3, see specifically Lemma~3.2,
we see that the
conditions of Theorem~2.5 are satisfied for algebras whose minimal projections
are not too big.
This may include ``most'' algebras of interest.
However, Theorem~2.5 is a far from satisfactory answer to the question
about factoriality.
Indeed, consider
$(\MvN,\psi)=(\Cpx\bigoplus\Cpx,\phi)*(M_2(\Cpx),\tau)$, where $\tau$ is
the trace--state on the $2\times2$ matrices and where $\phi$ is the state on
$\Cpx^2$ assigning to one minimal projection the weight $\lambda$,
$\frac12\le\lambda<1$.
{}From~\cite{D3} we have that $\MvN$ is a factor if and only if
$\frac12\le\lambda\le\frac34$.
However, we have calculated (Examples~3.4i and~3.4iii) that
$$ \ef(\Cpx^2,\phi)=\cases1&\text{if }\lambda=\frac12 \\
0&\text{if }\lambda>\frac12,
\endcases $$
and thus only in the case $\lambda=\frac12$ does
Theorem~2.5 imply that $\MvN$ is a factor.
\item"ii)"  Consider
the situation of Example~3.4iv, {\it i.e\.} $A=M_2(\Cpx)$ with
the state which we will denote $\phi_\lambda$.
Figure~1 is a plot of $\ef(A,\phi_\lambda)$ versus $\lambda$.
We thus find that
if $\lambda,\mu\in[\frac12,1)$ and if $(\MvN,\psi)=(A,\phi_\lambda)
*(A,\phi_\mu)$, then Theorem~2.5 implies that $\MvN$ is a factor for
$(\lambda,\mu)$ inside the region bounded be the curve in Figure~2.
One might speculate that $\MvN$ should be factor for all
choices of $\lambda$ and $\mu$,
and that Figure~2 illustrates
only the limitations of Theorem~2.5.
\endroster
\endproclaim

\proclaim{Examples 4.2} Type Classification. \rm
\roster
\item"i)" Let $(\MvN,\phi)$ and $(\NvN,\psi)$ be type~III von Neumann
algebras, each with separable predual.
Then their free product von Neumann
algebra is a type III factor.
\item"ii)" More generally, if $(A,\phi)$ and $(B,\psi)$ are any von Neumann
algebras having separable preduals and such that $\phi$ and $\psi$ are
normal faithful states, one of which is not a trace,
and if $A$ (respectively $B$) has no minimal projections $e$ such
that $\phi(e)>\frac12$ (respectively $\psi(e)>\frac12$),
then the free product von Neumann algebra $A*B$
of $(A*B,\phi*\psi)$ is a type III factor.
\item"iii)" If in the above case we assume that $\psi$ is a trace,
then $T(A*B)=I(A,\phi)$.
\item"iv)" If $(\MvN,\phi)$ and $(\NvN,\psi)$ are factors with separable
preduals and having types III$_\lambda$ and III$_\mu$, respectively,
with $0<\lambda,\mu\le1$, and if either $\lambda=1$ or $\mu=1$, or if
$\ln\lambda$ and $\ln\mu$ are not rationally related, then the free product
of $\MvN$ and $\NvN$ is a factor of type III$_1$ or type III$_0$.
If $\ln\lambda$ and
$\ln\mu$ are rationally related, then $\MvN$ is either type~III$_\nu$
or type~III$_0$, where $0<\nu<1$ is such that
$\frac{2\pi}{\ln\nu}\Integers=\frac{2\pi}{\ln\lambda}\Integers\cap
\frac{2\pi}{\ln\mu}\Integers$.
\endroster
\endproclaim
\demo{Proof}
Since, given two von Neumann algebras having separable preduals,
their free product von Neumann algebra is faithfully represented on
a separable Hilbert space, it too has separable predual.
We know from Theorem~2.5 and Lemma~3.2
that the free product von Neumann algebras
in the cases given above are factors, so their $T$--invariants
being not $\Real$ implies that
the factor is type III.
To prove~iii), note that the modular automorphism group defined with
respect to a trace is the identity, so $I(B,\psi)=\Real$.
In the cases described in part~iv, we have clearly that
the $T$-invariant of the free product factor equals $\{0\}$,
implying that the factor is type III$_1$ or III$_0$.
\QED

  If $\MvN$ has semifinite trace $\Tr$ and if $\phi$ is a normal
faithful state on $\MvN$ having density $h\in\MvN_+$, {\it i.e\.}
$\phi(x)=\Tr(hx)$, then~(\cite{T}, p.98), the modular automorphisms
group is given by $\sigma_t^\phi(x)=u(t)xu(-t)$ where
$u(t)=h^{it}$.
Thus $T(\MvN)=\Real$ and $I(\MvN,\phi)=\{t\in\Real\mid h^{it}\in\Cpx1\}$.
For the case of $n\times n$
matrix algebras, let $(A,\phi)$ be as described in Lemma~3.3,
and let $\Tr$ be the trace on $M_n(\Cpx)$ such that $\Tr(e_{\iota\iota})=1$.
Then the density $h$ is the diagonal matrix with diagonal entries
$\lambda_1,\lambda_2,\ldots,\lambda_n$.
Hence
$$ \align
I(A,\phi)&=\{t\in\Real\mid\left(\frac{\lambda_\iota}{\lambda_j}\right)^{it}
=1\;\forall\;1\le\iota,j\le n\} \\
&=\{t\in\Real\mid\left(\frac{\lambda_1}{\lambda_\iota}\right)^{it}
=1\;\forall\;1\le\iota\le n\}.
\endalign $$
So if $\phi$ is not a trace then
$$ I(A,\phi)=\bigcap_{\{2\le\iota\le n\mid\lambda_\iota\neq\lambda_1\}}
\left(\frac{2\pi}{\ln\lambda_1-\ln\lambda_\iota}\right)\Integers.
\tag23 $$

  To make a connection with F\. R\u{a}dulescu's work in~\cite{R5},
for the state $\psi_\lambda$ on $M_2(\Cpx)$ (described in the introduction),
we get from~(23) that
$I(M_2(\Cpx),\psi_\lambda)=\frac{2\pi}{-\ln\lambda}\Integers$.
Then Theorem~2.5 implies that for
$(\MvN,\phi)=(M_2(\Cpx),\psi_\lambda)*(L(\Integers),\tau)$, where $\tau$
is a faithful trace on the diffuse Abelian von Neumann algebra $L(\Integers)$,
we have that $\MvN$ is a type~III factor and that
$T(\MvN)=\frac{2\pi}{-\ln\lambda}\Integers$, so that $\MvN$ is
either type~III$_\lambda$ or type~III$_0$.
R\u{a}dulescu has shown that it is type~III$_\lambda$ and has further described
$\MvN$ in terms of a discrete decomposition for it.

  When the conditions of Theorem~2.5 are satisfied for a free product
of algebras with states for which we can compute the $I$--invariants,
we can use~(23) and~(7) to make conclusions about the $\lambda$--classification
of the free product factor,
similarly to what we did in Example~4.2iv.  For example,
in the notation of the previous paragraph, if
$(\MvN,\phi)=(M_2(\Cpx),\psi_\lambda)*(M_s(\Cpx),\psi_\mu)$,
then we can as in Remark~4.1ii decide when the conditions of Theorem~2.5
are satisfied, (note that $\psi_\lambda=\phi_{\lambda'}$,
where $\lambda'=\frac1{1+\lambda}$).
Suppose the conditions of Theorem~2.5 are satisfied.
If $\ln\lambda$ and
$\ln\mu$ are not rationally related, then $\MvN$ is either type~III$_1$
or type~III$_0$.
If $\ln\lambda$ and
$\ln\mu$ are rationally related, then $\MvN$ is either type~III$_\nu$
or type~III$_0$, where $0<\nu<1$ is such that
$\frac{2\pi}{\ln\nu}\Integers=\frac{2\pi}{\ln\lambda}\Integers\cap
\frac{2\pi}{\ln\mu}\Integers$.

\end{document}